\documentclass[12pt,twoside]{article}
\usepackage{amsmath,amssymb}
\usepackage{graphicx}
\usepackage{epstopdf}
\usepackage{float}
\usepackage{mathrsfs,amssymb}

\allowdisplaybreaks[4]

\newtheorem{lm}{Lemma}[section]
\newtheorem{thm}{Theorem}[section]
\newtheorem{rmk}{Remark}[section]

\newtheorem{exa}{Example}[section]

\newcounter{saveeqn}%

\makeatletter
\evensidemargin
\oddsidemargin
\makeatother

\title{\Large\bf
Conjugacy problem of strictly monotone maps with only one jump discontinuity
\thanks{The first author is supported by the general item of Lingnan Normal University [grant ZL1505], KSP of Lingnan Normal University [No.1171518004].
The second author is supported by Scientific Research Fund of SiChuan Provincial Education Department (18ZA0274) and NNSF (11301256). }}
\author{Jinghua Liu $^{a}$,~~
Yong-Guo Shi $^{b}$
\footnote{Corresponding author.
E-mail addresses: scumat@163.com (Y. G. Shi).
}
\\
$^{a}${\small School of Mathematics and Statistics, Lingnan Normal University }
\\
{\small Zhanjiang, Guangdong 524048, P.R.China}
\\
$^{b}${\small  Data Recovery Key Laboratory of Sichuan Province,
}
\\
{\small College of Mathematics and Information Science, Neijiang Normal University,}
\\
{\small Neijiang, Sichuan 641112, P.R.China}}

\date{ }

\begin{document}
\maketitle

\begin{abstract}

The conjugacy problem is one of the central questions in iteration theory.
As far as we, for discontinuous strictly monotone maps there is no complete result.
In this paper, we investigate the conjugacy problem of strictly monotone maps
with only one jump discontinuity. We give some sufficient and necessary conditions for the conjugacy relationship.
And we present some methods to construct all conjugacies. Furthermore, we present the conditions to guarantee $C^1$ smoothness of these conjugacies.
\vskip 0.2cm
{\bf Keywords:} discontinuous interval map; conjugacy; jump discontinuity; smoothness.
\vskip 0.2cm
{\bf AMS(2010) Subject Classifications:} 39C15, 37E05
\end{abstract}

\baselineskip 15pt
\parskip 10pt
\thispagestyle{empty}
\setcounter{page}{1}


\section{Introduction}
\setcounter{equation}{0}
\setcounter{lm}{0}
\setcounter{thm}{0}
\setcounter{rmk}{0}
\setcounter{df}{0}
\setcounter{cor}{0}
\setcounter{exa}{0}

Let $I$, $J$ be two closed intervals, and $f : I\to I$ and $g : J \to J$ be two self-maps.
We say that $f$ and $g$ are topologically conjugate if there exists a homeomorphism $\varphi: I \to J$ satisfying
the conjugacy equation
\begin{eqnarray}
\varphi\circ f=g\circ\varphi.
\label{conjugacy equ}
\end{eqnarray}
Here $\varphi$ is called a conjugacy from $f$ to $g$.

The conjugacy relation is an important equivalent relation in dynamical systems,
and it is used for topological classification and simplification of dynamical systems and functional equations.
In general, the conjugacy problem contains three aspects. The first one is how to determine any two maps in some family
to be topologically conjugate.
The second is how to construct all conjugacies. The last is what about the smooth of every conjugacy.
Many works are devoted to the conjugacy problem of continuous interval maps (i.e., \cite{Block-Coven1987,Jiang1995,Lesniak-Shi2015,Li-Shen2006,Parry1966,Segawa-Ishitani1998,Shi2009}).

However, there is a few work for the conjugacy problem of discontinuous maps, for example, Lorenz map, cf. \cite{Cui-Ding2015,Glendinning1990,Glendinning-Sparrow1993,Hubbard-Sparrow1990,Llibre1987,Pring-Budd2010}.
As far as we, even for discontinuous strictly monotone maps there is no complete result.
Consider a family of simple discontinuous strictly monotone maps
with only one jump discontinuity. It is interesting to consider how this jump discontinuity affects the conjuacy relationship.
One can see that there exist three main factors: (i) the value of functions at the jump discontinuity, (ii) the position relationship
between the image of functions and the diagonal, (iii) the right and left limits at the discontinuity.

In this paper, we firstly give some necessary conditions for the cojugacy relationship in the next section. Further, we distinguish between increasing conjugacies and decreasing conjugacies. Section 3 is devoted for sufficient and necessary conditions for these discontinuous maps, and construct all conjugacies.
In section 4, we give some sufficient conditions under which these conjugacies are $C^{1}$ smooth. In the final section, two examples are presented.

\section{Preliminaries}
\setcounter{equation}{0}
\setcounter{lm}{0}
\setcounter{thm}{0}
\setcounter{rmk}{0}
\setcounter{df}{0}
\setcounter{cor}{0}
\setcounter{exa}{0}

We only consider the following representative cases, others can be discussed similarly.
Let $I:=[a,b]$ and $\mathcal{A}_{t}(I)$ denote the family of strictly increasing functions $f$ with only one
jump discontinuous point $t\in (a,b)$ satisfying (i) $f(t-0)\leq f(t)\leq f(t+0)$ and (ii) $f$ has exactly two attractive fixpoints $a$ and $b$.
See Figs 1-4.

\begin{figure}[H]
\begin{minipage}[t]{0.5\linewidth}
\centering
\includegraphics[width=0.7\textwidth]{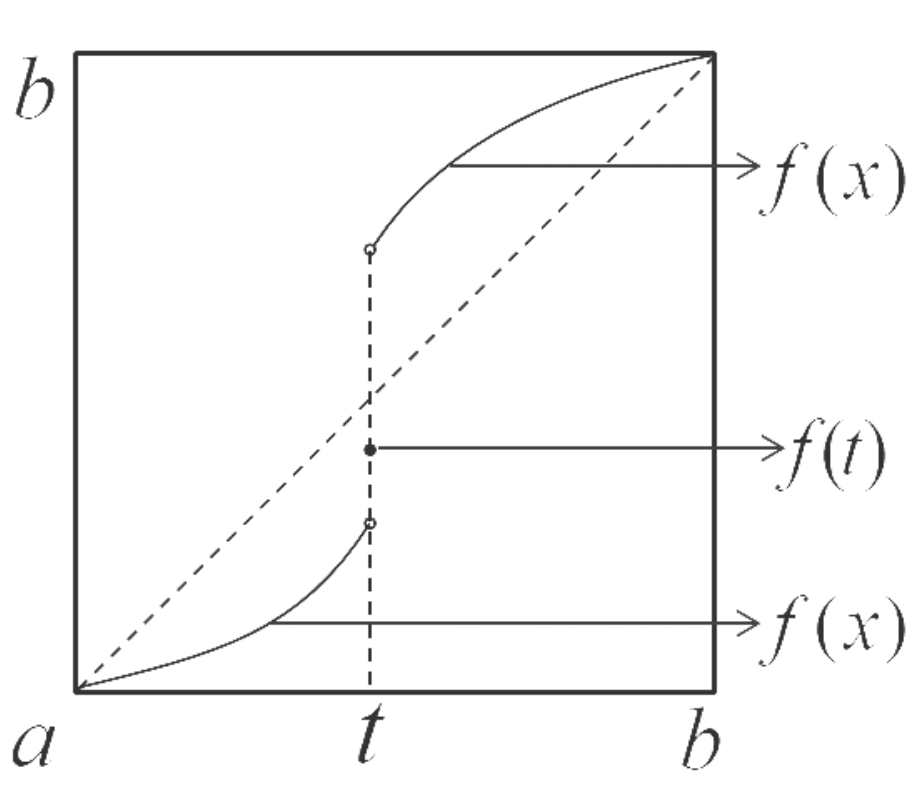}
\caption{$f\in\mathcal{A}_{t}(I),~f(t)<t$} \label{fig 1}
\end{minipage}%
\begin{minipage}[t]{0.5\linewidth}
\centering
\includegraphics[width=0.7\textwidth]{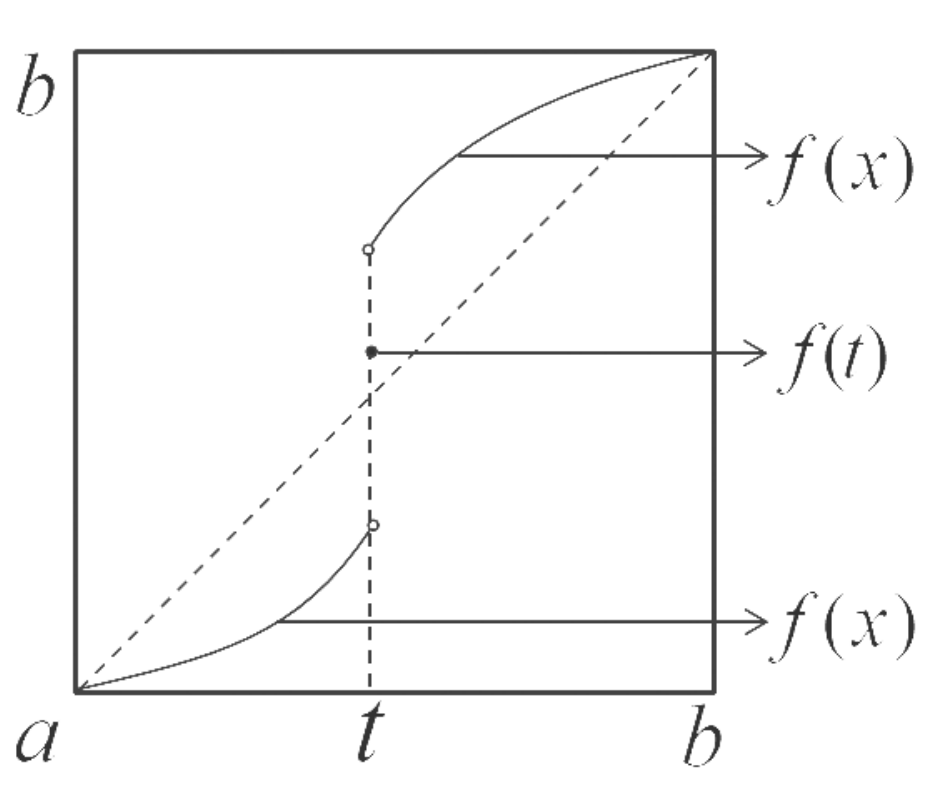}
\caption{$f\in\mathcal{A}_{t}(I),~f(t)>t$} \label{fig 2}
\end{minipage}

\begin{minipage}[t]{0.5\linewidth}
\centering
\includegraphics[width=0.7\textwidth]{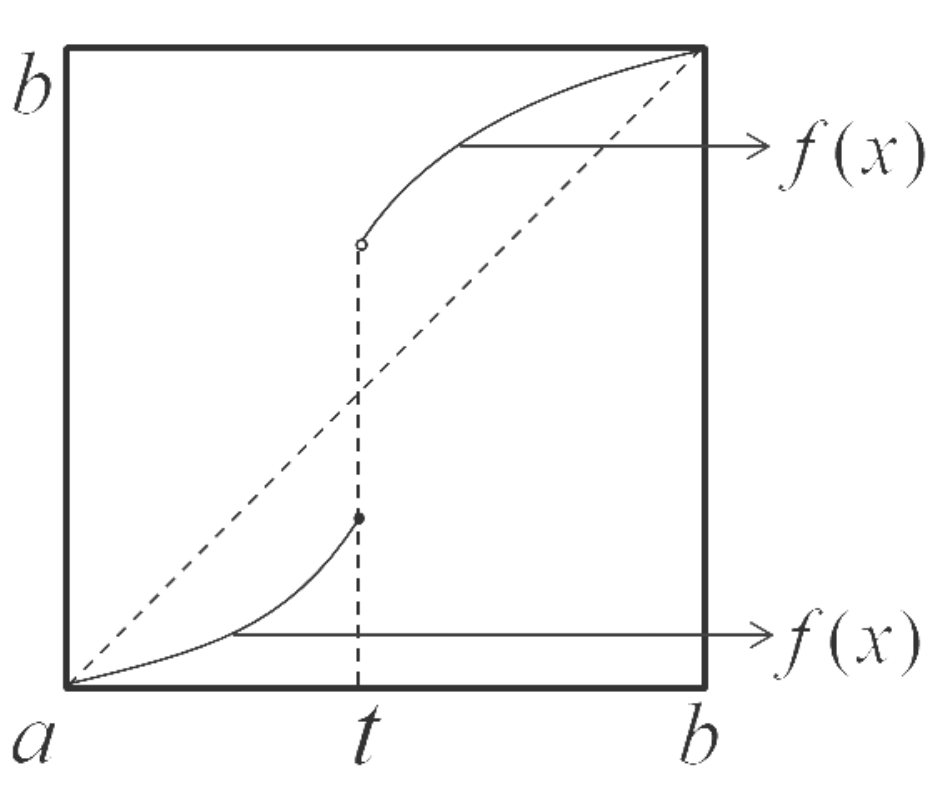}
\caption{$f\in\mathcal{A}_{t}(I),~f(t)<t$} \label{fig 3}
\end{minipage}%
\begin{minipage}[t]{0.5\linewidth}
\centering
\includegraphics[width=0.7\textwidth]{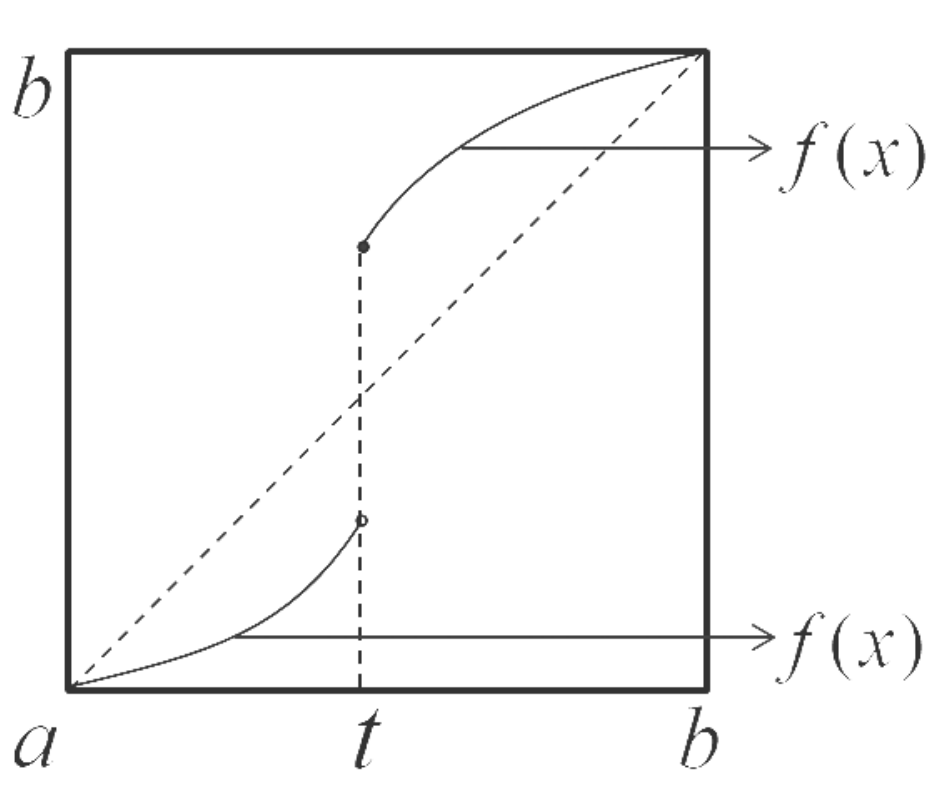}
\caption{$f\in\mathcal{A}_{t}(I),~f(t)>t$} \label{fig 4}
\end{minipage}
\end{figure}
Let $\mathcal{B}_{t}(I)$ denote the family of strictly decreasing functions $f$ with only one
jump discontinuous point $t\in (a,b)$ satisfying (i) $f(t-0)\leq f(t)\leq f(t+0)$ and (ii) $f$ has exactly two periodic points $a$ and $b$,
both of which are period~2~points. See Figs 5-8.

\begin{figure}[H]
\begin{minipage}[t]{0.5\linewidth}
\centering
\includegraphics[width=0.7\textwidth]{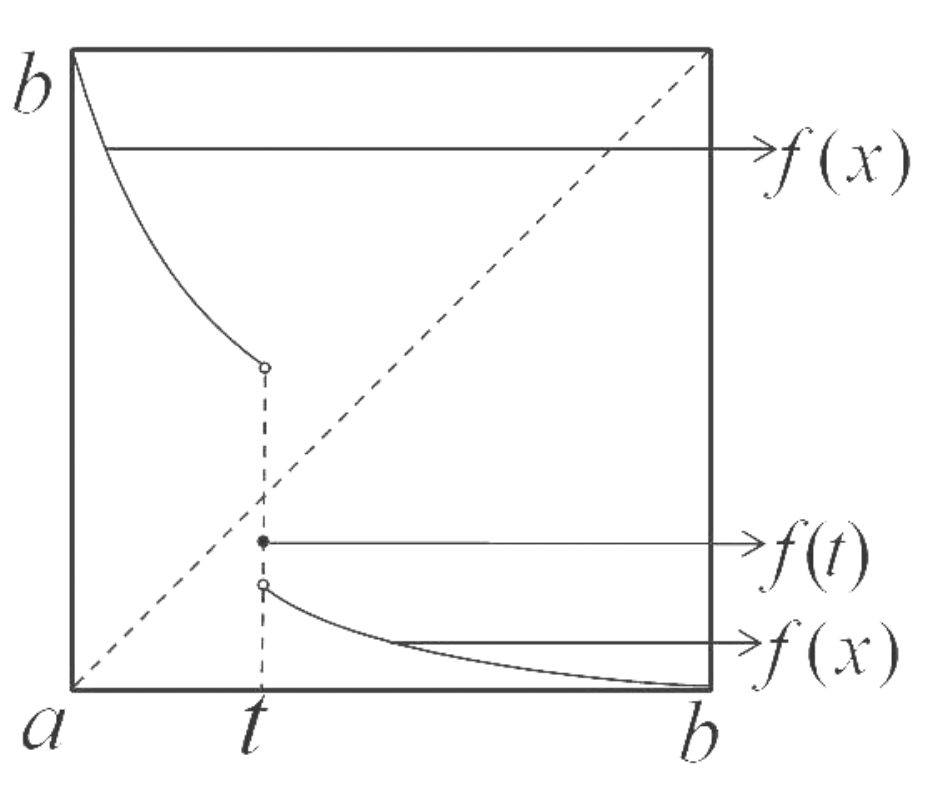}
\caption{$f\in\mathcal{B}_{t}(I),~f(t)<t$} \label{fig 5}
\end{minipage}%
\begin{minipage}[t]{0.5\linewidth}
\centering
\includegraphics[width=0.7\textwidth]{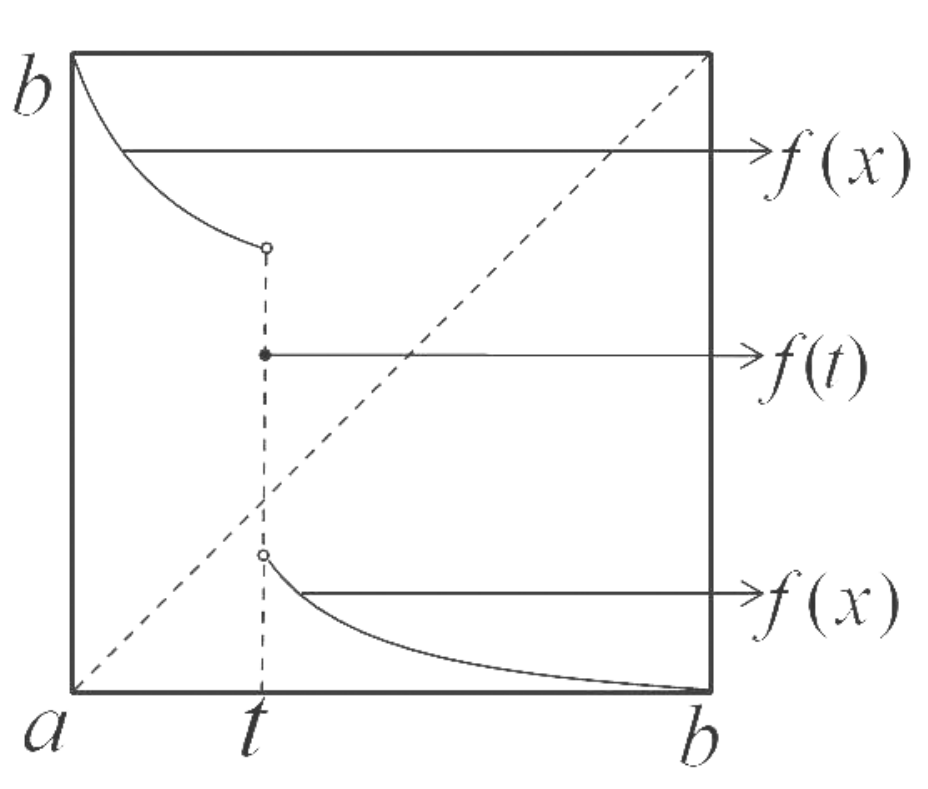}
\caption{$f\in\mathcal{B}_{t}(I),~f(t)>t$} \label{fig 6}
\end{minipage}
\end{figure}
\begin{figure}[H]
\begin{minipage}[t]{0.5\linewidth}
\centering
\includegraphics[width=0.7\textwidth]{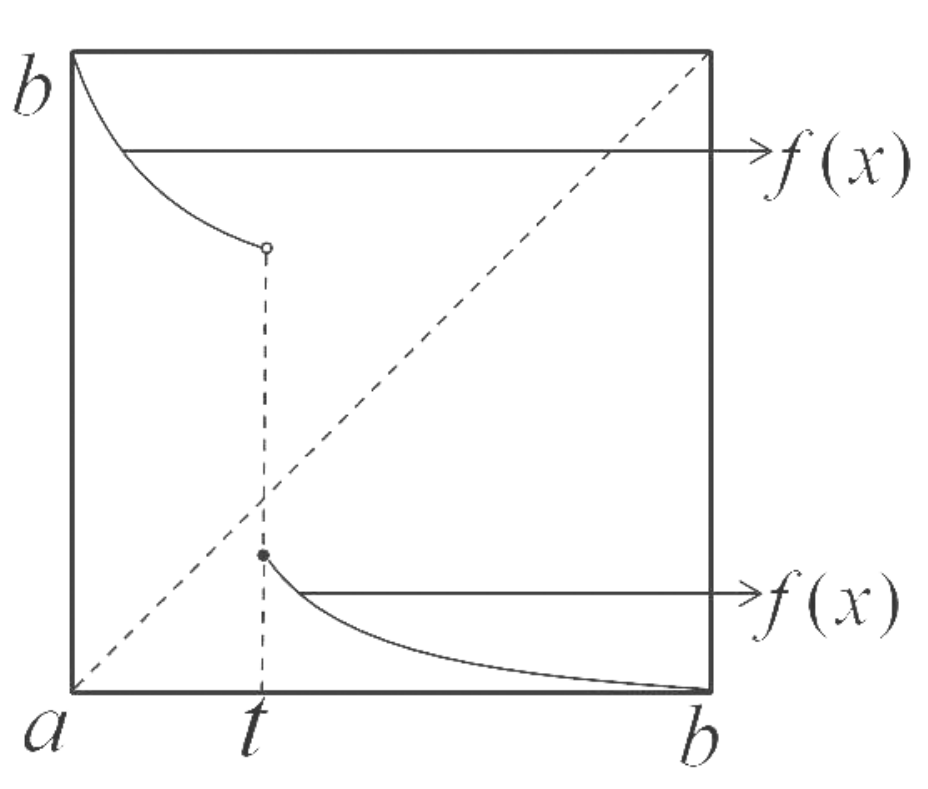}
\caption{$f\in\mathcal{B}_{t}(I),~f(t)<t$} \label{fig 7}
\end{minipage}%
\begin{minipage}[t]{0.5\linewidth}
\centering
\includegraphics[width=0.7\textwidth]{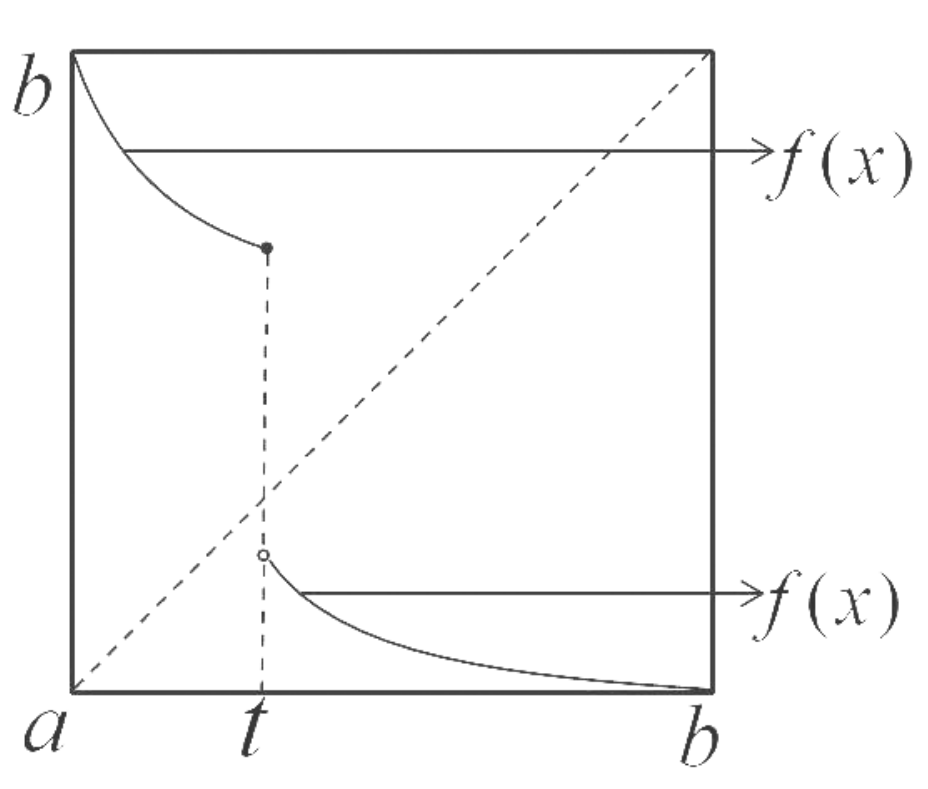}
\caption{$f\in\mathcal{B}_{t}(I),~f(t)>t$} \label{fig 8}
\end{minipage}
\end{figure}

Let $J:=[c,d]$ and $c<s<d$. Some necessary conditions for the conjugacy relationship are given below.
\begin{lm}
Let $f\in\mathcal{A}_{t}(I)$ and $g\in\mathcal{A}_{s}(J)$ be topologically conjugate via a homeomorphism $\varphi$.
Then $\varphi(f^n(t))=g^n(s)$ for $n=0,1,2,...$.
\label{lemma4}
\end{lm}
\noindent{\bf Proof.} Since $\varphi\circ f$ is not continuous at point $t$, and $g\circ \varphi=\varphi\circ f$, $g$ is not continuous at $\varphi(t)$. Therefore $\varphi(t)=s$. Further, $\varphi(f(t))=g(\varphi(t))=g(s)$. By induction, $\varphi(f^n(t))=g^n(s)$ for $n=0,1,2,...$. $\Box$

\begin{lm}
Let $f\in\mathcal{A}_{t}(I)$ and $g\in\mathcal{A}_{s}(J)$ be topologically conjugate. Then
there are only two cases:
\begin{description}
\item[(C1)]
$f(t-0)<f(t)<f(t+0)$ and $g(s-0)<g(s)<g(s+0)$;
\item[(C2)]
$f(t)=f(t-0)$ and $g(s)=g(s-0)$, or $f(t)=f(t-0)$ and $g(s)=g(s+0)$, or
$f(t)=f(t+0)$ and $g(s)=g(s-0)$, or $f(t)=f(t+0)$ and $g(s)=g(s+0)$.
\end{description}
\label{lemma3}
\end{lm}
\noindent{\bf Proof.}
We only prove the case $f(t-0)< f(t)<f(t+0)$ and $g(s)=g(s-0)$. Others are similar.
 Without loss of generality, assume that $\varphi$ is an increasing conjugacy from $f$ to $g$.
One can see that $\varphi(f(t-0))=g(\varphi(t-0))=g(s)=g(s-0)$.
 It follows from Lemma \ref{lemma4} that $\varphi(f(t))=g(\varphi(t))=g(s)=g(s-0)$. Then
$\varphi(f(t))=\varphi(f(t-0))$. Since $\varphi$ is a  homeomorphism, we have $f(t-0)=f(t)$. This is a contradiction.
 $\Box$

\begin{lm}
Let $f\in\mathcal{A}_{t}(I)$ and $g\in\mathcal{A}_{s}(J)$ be topologically conjugate via a homeomorphism $\varphi$.
\begin{description}
\item[(i)]
If $(f(t)-t)(g(s)-s)>0$, then $\varphi$ is increasing and goes through points $(a,c)$, $(b,d)$, $(t,s)$, $(f(t),g(s))$,  $(f(t-0),g(s-0))$ and $(f(t+0),g(s+0))$.

\item[(ii)]
If $(f(t)-t)(g(s)-s)<0$, then $\varphi$ is decreasing and goes through points $(a,d)$, $(b,c)$, $(t,s)$, $(f(t),g(s))$,  $(f(t+0),g(s-0))$ and $(f(t-0),g(s+0))$.

\item[(iii)]
If $f(t)=f(t\pm 0)$ and $g(t)=g(t\pm 0)$,  then $\varphi$ is increasing and goes through points $(a,c)$, $(b,d)$, $(t,s)$, $(f(t-0),g(s-0))$ and $(f(t+0),g(s+0))$.

\item[(iv)]
If $f(t)=f(t\pm 0)$ and $g(t)=g(t\mp 0)$, then $\varphi$ is decreasing and goes through points $(a,d)$, $(b,c)$, $(t,s)$, $(f(t+0),g(s-0))$ and $(f(t-0),g(s+0))$.
\end{description}
\label{lemma1}
\end{lm}
\noindent{\bf Proof.}
We only prove the fact \textbf{(i)}. Others are similar.
For fact \textbf{(i)},  we will prove that $\varphi$ is increasing by contradiction.
 We only consider the case $f(t)<t$ and $g(s)<s$, the other case is similar.
Assume that $\varphi$ is decreasing, then $\varphi(f(t))>\varphi(t)$. Thus
$g(\varphi(t))>\varphi(t)$, i.e., $g(s)>s$, which contradicts the condition that $g(s)<s$. Therefore, $\varphi$ is increasing.

From $\varphi(f(a))=g(\varphi(a))$, we have $\varphi(a)=g(\varphi(a))$. Thus $\varphi(a)=c$ or $\varphi(a)=d$. Since $\varphi$ is strictly increasing, it is clear that $\varphi(a)=c$. Similarly, we have $\varphi(b)=d$.

Since $\lim_{x\rightarrow t^{-}}\varphi(f(x))=\lim_{x\rightarrow t^{-}}g(\varphi(x))$, we have
$$
\varphi(\lim_{x\rightarrow t^{-}}f(x))=\lim_{x\rightarrow t^{-}}g(\varphi(x)).
$$
Then $\varphi(f(t-0))=g(s-0)$. Similarly, we get $\varphi(f(t+0))=g(s+0)$.
$\Box$

For the family $\mathcal{B}_{t}$, we can obtain the corresponding results by the similar arguments.
\begin{lm}
Let $f\in\mathcal{B}_{t}(I)$ and $g\in\mathcal{B}_{s}(J)$ be topologically conjugate via a homeomorphism $\varphi$.
Then $\varphi(t)=s$ and $\varphi(f(t))=g(s)$.
\label{lemma6}
\end{lm}
\begin{lm}
Let $f\in\mathcal{B}_{t}(I)$ and $g\in\mathcal{B}_{s}(J)$ be topologically conjugate. Then
there are only two cases:
\begin{description}
\item[(D1)]
$f(t-0)<f(t)<f(t+0)$ and $g(s-0)<g(s)<g(s+0)$;
\item[(D2)]
$f(t)=f(t-0)$ and $g(s)=g(s-0)$, or $f(t)=f(t-0)$ and $g(s)=g(s+0)$, or
$f(t)=f(t+0)$ and $g(s)=g(s-0)$, or $f(t)=f(t+0)$ and $g(s)=g(s+0)$.
\end{description}
\label{lemma5}
\end{lm}

\begin{lm}
Let $f\in\mathcal{B}_{t}(I)$ and $g\in\mathcal{B}_{s}(J)$ be topologically conjugate via a homeomorphism $\varphi$.
\begin{description}
\item[(i)]
If  $(f(t)-t)(g(s)-s)>0$, then $\varphi$ is increasing and goes through points $(a,c)$, $(b,d)$, $(t,s)$, $(f(t),g(s))$,
$(f^{2}(t-0),g^{2}(s-0))$, $(f^{2}(t+0),g^{2}(s+0))$, $(f(t+0),g(s+0))$, and $(f(t-0),g(s-0))$.

\item[(ii)]
If  $(f(t)-t)(g(s)-s)<0$, then $\varphi$ is decreasing and goes through points $(a,d)$, $(b,c)$, $(t,s)$, $(f(t),g(s))$,  $(f^{2}(t+0),g^{2}(s-0))$,
$(f^{2}(t-0),g^{2}(s+0))$, $(f(t-0),g(s+0))$, and $(f(t+0),g(s-0))$.

\item[(iii)]
If $f(t)=f(t\pm 0)$ and $g(t)=g(t\pm 0)$,  then $\varphi$ is increasing and goes through points $(a,c)$, $(b,d)$, $(t,s)$,
$(f^{2}(t-0),g^{2}(s-0))$, $(f^{2}(t+0),g^{2}(s+0))$, $(f(t+0),g(s+0))$, and $(f(t-0),g(s-0))$.

\item[(iv)]
If $f(t)=f(t\pm 0)$ and $g(t)=g(t\mp 0)$,  then $\varphi$ is decreasing and goes through points $(a,d)$, $(b,c)$, $(t,s)$, $(f^{2}(t+0),g^{2}(s-0))$,
$(f^{2}(t-0),g^{2}(s+0))$, $(f(t-0),g(s+0))$, and $(f(t+0),g(s-0))$.
\end{description}
\label{lemma2}
\end{lm}

\section{Construction of topological conjugacy}
\setcounter{equation}{0}
\setcounter{lm}{0}
\setcounter{thm}{0}
\setcounter{rmk}{0}
\setcounter{df}{0}
\setcounter{cor}{0}
\setcounter{exa}{0}

In this section, we give the sufficient and necessary conditions for the conjugacy in the family $\mathcal{A}_{t}$ and $\mathcal{B}_{t}$, and construct all conjugacies.

\subsection{The family $\mathcal{A}_{t}$}

\begin{thm}
Suppose that $f\in\mathcal{A}_{t}(I)$ and $g\in\mathcal{A}_{s}(J)$. Then  $f$ and $g$ are topologically conjugate if and only if $f$ and $g$ are in case ({\bf C$1$}) or both in ({\bf C$2$}).
\label{increasing-case-l}
\end{thm}

\noindent{\bf Proof.}
The necessarity follows from Lemma \ref{lemma3}. It suffices to prove the sufficiency.

We only consider case \textbf{(C1)}. The other is similar. For case  \textbf{(C1)}, we have the following four subcases:
\begin{eqnarray*}
&&\textbf{(i)}~f(t)<t~{\rm and}~g(s)<s,~\textbf{(ii)}~f(t)<t~{\rm and}~g(s)>s,\\
&&\textbf{(iii)}~f(t)>t~{\rm and}~g(s)<s,~\textbf{(iv)}~f(t)>t~{\rm and}~g(s)>s.
\end{eqnarray*}
We only consider subcases \textbf{(i)}-\textbf{(ii)}. Other subcases can be proved similarly.

For subcase \textbf{(i)}.
Firstly, define functions $f_{l}$, $f_{r}$, $g_{l}$ and $g_{r}$ as follows
\begin{eqnarray*}
f_{l}(x):\!\!\!\! &=& \!\!\!\! \left\{
\begin{array}{ll}
f(x),           &x\in [a, t),
\\
f(t-0), & x=t,
\end{array}
\right.
\\
f_{r}(x):\!\!\!\! &=& \!\!\!\! \left\{
\begin{array}{ll}
f(x),           &x\in (t,b],
\\
f(t+0),         & x=t,
\end{array}
\right.
\\
g_{l}(x):\!\!\!\! &=& \!\!\!\! \left\{
\begin{array}{ll}
g(x),           &x\in [c,s),
\\
g(s-0), & x=s,
\end{array}
\right.
\\
g_{r}(x):\!\!\!\! &=& \!\!\!\! \left\{
\begin{array}{ll}
g(x),           &x\in (s,d],
\\
g(s+0), & x=s.
\end{array}
\right.
\end{eqnarray*}
Then define two strictly monotone sequences $\{f_{l}^{n}(t)\}$ and $\{f_{r}^{n}(t)\}$, where $n=0,1,2,\cdots$. Obviously, $f_{l}^{n}(t)\rightarrow a$ and $f_{r}^{n}(t)\rightarrow b$ as $n\rightarrow +\infty$. With these sequences, give a partition of $I$, i.e.,
$$(a,~t]=\cup_{n=0}^{+\infty}I_{n},~~[t,~b)=\cup_{n=0}^{+\infty}J_{n},$$
where $I_{n}=(f_{l}^{n+1}(t),f_{l}^{n}(t)]$ and $J_{n}=[f_{r}^{n}(t),f_{r}^{n+1}(t))$.
By means of Lemma \ref{lemma1}, arbitrarily choose two increasing homeomorphisms $\varphi_{l_{0}}:[f_{l}(t),t]\rightarrow [g_{l}(s),s] $ and $\varphi_{r_{0}}:[t,f_{r}(t)]\rightarrow [s,g_{r}(s)] $ as initial functions such that $\varphi_{l_{0}}(f(t))=g(s)$. With the iterative constructive method and Lemma \ref{lemma1}, define for $n=1,2,...$
\begin{eqnarray}
\varphi_{l}(x):\!\!\!\! &=& \!\!\!\! \left\{
\begin{array}{ll}
\varphi_{l_{0}}(x),           &x\in (f_{l}(t),  t],
\\
g_{l}^{n}\circ\varphi_{l_{0}}\circ f_{l}^{-n}(x), &x\in (f_{l}^{n+1}(t),  f_{l}^{n}(t)],
\\
c, &x=a,
\end{array}
\right.
\label{inc-solu1}
\end{eqnarray}
\begin{eqnarray}
\varphi_{r}(x):\!\!\!\! &=& \!\!\!\! \left\{
\begin{array}{ll}
\varphi_{r_{0}}(x),           &x\in [t,  f_{r}(t)),
\\
g_{r}^{n}\circ\varphi_{r_{0}}\circ f_{r}^{-n}(x), &x\in [f_{r}^{n}(t), f_{r}^{n+1}(t)),
\\
d, &x=b.
\end{array}
\right.
\label{inc-solu2}
\end{eqnarray}
Thus $\varphi_{l}: [a,t]\to [c,s]$ is an increasing conjugacy from $f_{l}$  to $g_{l}$ and $\varphi_{r}: [t,b]\to [s,d]$ is an increasing conjugacy from $f_{r}$ to $g_{r}$.
In fact, for $x\in I_{n}$,
it follows that
\begin{eqnarray}
\varphi_{l}(f_{l}(x))=g_{l}^{n+1}\circ \varphi_{l_{0}}\circ f_{l}^{-n-1}(f_{l}(x))=
g_{l}\circ g_{l}^{n}\circ \varphi_{l_{0}}\circ f_{l}^{-n}(x)=g_{l}(\varphi_{l}(x)),
\label{conjugcy-equ-trans}
\end{eqnarray}
\begin{eqnarray}
\nonumber
\varphi_{l}(f_{l}(a))=g(\varphi_{l}(a)).
\end{eqnarray}
Now, it suffices to show $\varphi_{l}$ is strictly increasing and continuous on the interval $[a,t]$.
One can see $\varphi_{l}$ is strictly increasing and continuous on each open interval $(f_{l}^{n+1}(t),  f_{l}^{n}(t))$ for $n\in\mathbb{N}$.
Thus it suffices to show $\varphi_{l}$ is continuous at each transitional point $f_{l}^{n}(t)$ for $n\in\mathbb{N}$ and right-continuous at $x=a$. By induction, for $n=1$, we have
$$
\lim_{h\rightarrow 0^{+}}\varphi_{l}(f_{l}(t)+h)=\lim_{h\rightarrow 0^{+}}\varphi_{l_{0}}(f_{l}(t)+h)=\varphi_{l_{0}}(f_{l}(t))=\varphi_{l}(f_{l}(t)),
$$
$$
\lim_{h\rightarrow 0^{+}}\varphi_{l}(f_{l}(t)-h)=\lim_{h\rightarrow 0^{+}}g_{l}\circ\varphi_{l_{0}}\circ f_{l}^{-1}(f_{l}(t)-h)=g_{l}\circ\varphi_{l_{0}}\circ f_{l}^{-1}(f_{l}(t))=\varphi_{l}(f_{l}(t)).
$$
Then assume $\varphi_{l}$ is continuous at the point $f_{l}^{n}(t)$ for some positive integer $n$.
We shall show $\varphi_{l}$ is continuous at the transitional point $f_{l}^{n+1}(t)$.
For convenience, we denote $g_{l}^{n}\circ\varphi_{l_{0}}\circ f_{l}^{-n}$ by $\varphi_{l_{n}}$. Thus $\varphi_{l}=\varphi_{l_{n}}$ on the interval $(f_{l}^{n+1}(t),~f_{l}^{n}(t)]$.
Since
$$
\lim_{h\rightarrow 0^{+}}\varphi_{l}(f_{l}^{n+1}(t)+h)=\lim_{h\rightarrow 0^{+}}\varphi_{l_{n}}(f_{l}^{n+1}(t)+h)=\varphi_{l_{n}}(f_{l}^{n+1}(t))=\varphi_{l}(f_{l}^{n+1}(t))
$$
and
\begin{eqnarray*}
\lim_{h\rightarrow 0^{+}}\varphi_{l}(f_{l}^{n+1}(t)-h)&=&\lim_{h\rightarrow 0^{+}}g_{l}\circ\varphi_{l_{n}}\circ f_{l}^{-1}(f_{l}^{n+1}(t)-h)\\
&=&g_{l}\circ\varphi_{l_{n}}\circ f_{l}^{-1}(f_{l}^{n+1}(t))\\
&=& \varphi_{l}(f_{l}^{n+1}(t)),
\end{eqnarray*}
Thus $\varphi_{l}$ is continuous at the transitional point $f_{l}^{n+1}(t)$.
Therefore, $\varphi_{l}$ is continuous at the point $f_{l}^{n}(t)$ for $n\in\mathbb{N}$.
Finally, we shall show that $\varphi_{l}$ is right-continuous at $x=a$. Since $\varphi_{l}(t)=\varphi_{l_{0}}(t)=s$ and $\lim_{n\rightarrow\infty} g_{l}^{n}(s)=c$, we have $\varphi_{l}(f_{l}^{n}(t))=g_{l}^{n}(\varphi(t))$ and
\begin{align}
\lim_{n\rightarrow\infty}\varphi_{l}(f_{l}^{n}(t))=\lim_{n\rightarrow\infty}g_{l}^{n}(\varphi_{l}(t))=\lim_{n\rightarrow\infty}g_{l}^{n}(s)=c=\varphi(a).
\label{LD1}
\end{align}
Thus that $\varphi_{l}$ is right-continuous at $x=a$. Therefore $\varphi_{l}$ is continuous on $[a,t]$.

With the similar argument, we can prove that $\varphi_{r}: [t,b]\to [s,d]$ is an increasing conjugacy from $f_{r}$ to $g_{r}$.

Then define $\varphi$ on $[a,b]$ by
\begin{eqnarray}
\varphi(x):\!\!\!\! &=& \!\!\!\! \left\{
\begin{array}{ll}
\varphi_{l}(x),         &x\in [a,  t),
\\
s,                      &x=t,
\\
\varphi_{r}(x),         &x\in (t,  b],
\end{array}
\right.
\label{incrconj}
\end{eqnarray}
It is clear that $\varphi:I \to J$ is an increasing conjugacy from $f$ to $g$.

Finally, we will show all conjugacies can be obtained in this manner.

Suppose $\varphi$ is an increasing conjugacy between $f$ and $g$. Putting $\varphi_{l_{n}}:=\varphi$ on $I_{n}$ and $\varphi_{r_{n}}:=\varphi$ on $J_{n}$. By Lemma \ref{lemma1}, we can verify $\varphi_{l_{0}}$ goes through point $(t,s)$, $(f(t),g(s))$ and $(f_{l}(t),g_{l}(s))$, $\varphi_{r_{0}}$ goes through point $(t,s)$ and  $(f_{r}(t),g_{r}(s))$, and $\varphi_{l_{n}}$ (resp. $\varphi_{r_{n}}$) is of form (\ref{inc-solu1}) (resp. (\ref{inc-solu2})). Thus relation (\ref{incrconj}) holds.

For subcase \textbf{(ii)}, similarly, but choose any two decreasing homeomorphisms $\tilde{\varphi}_{l_{0}}:[f_{l}(t),t]\rightarrow [s,g_{r}(s)] $ and $\tilde{\varphi}_{r_{0}}:[t,f_{r}(t)]\rightarrow [g_{l}(s),s] $ as initial functions such that $\tilde{\varphi}_{l_{0}}(f(t))=g(s)$. Define for $n=1,~2,~3, \cdots$,
\begin{eqnarray*}
\tilde{\varphi}_{l}(x):\!\!\!\! &=& \!\!\!\! \left\{
\begin{array}{ll}
\tilde{\varphi}_{l_{0}}(x),           &x\in (f_{l}(t),  t],
\\
g_{r}^{n}\circ\tilde{\varphi}_{l_{0}}\circ f_{l}^{-n}(x), &x\in (f_{l}^{n+1}(t),  f_{l}^{n}(t)],
\\
d, &x=a,
\end{array}
\right.
\end{eqnarray*}
\begin{eqnarray*}
\tilde{\varphi}_{r}(x):\!\!\!\! &=& \!\!\!\! \left\{
\begin{array}{ll}
\tilde{\varphi}_{r_{0}}(x),           &x\in [t,  f_{r}(t)),
\\
g_{l}^{n}\circ\tilde{\varphi}_{r_{0}}\circ f_{r}^{-n}(x), &x\in [f_{r}^{n}(t), f_{r}^{n+1}(t)),
\\
c, &x=b,
\end{array}
\right.
\end{eqnarray*}
Then define
\begin{eqnarray*}
\tilde{\varphi}(x):\!\!\!\! &=& \!\!\!\! \left\{
\begin{array}{ll}
\tilde{\varphi}_{l}(x),         &x\in [a,  t),
\\
s,                      &x=t,
\\
\tilde{\varphi}_{r}(x),         &x\in (t,  b],
\end{array}
\right.
\end{eqnarray*}
which is a conjugacy from $f$ and $g$.
$\Box$

\subsection{The family $\mathcal{B}_{t}$}

\begin{thm}
Suppose that $f\in\mathcal{B}_{t}(I)$ and $g\in\mathcal{B}_{s}(J)$. Then  $f$ and $g$ are topologically conjugate if and only if $f$ and $g$ are in case ({\bf D$1$}) or both in ({\bf D$2$}).
\label{decreasing-case-1}
\end{thm}
\noindent{\bf Proof.}
The necessarity follows from Lemma \ref{lemma5}. It suffices to prove the sufficiency.

We only consider case \textbf{(D1)}. The other is similar. For case \textbf{(D1)}, we have the following four subcases:
\begin{eqnarray*}
&&\textbf{(i)}~f(t)<t~{\rm and}~g(s)<s,~\textbf{(ii)}~f(t)<t~{\rm and}~g(s)>s,\\
&&\textbf{(iii)}~f(t)>t~{\rm and}~g(s)<s,~\textbf{(iv)}~f(t)>t~{\rm and}~g(s)>s.
\end{eqnarray*}
We also prove subcase \textbf{(i)}-\textbf{(ii)}. Other subcases can be proved similarly.

For subcase \textbf{(i)}, define functions $f_{e}$, $f_{v}$, $g_{e}$ and $g_{v}$ as follows
\begin{eqnarray*}
f_{e}(x):\!\!\!\! &=& \!\!\!\! \left\{
\begin{array}{ll}
f(x),           &x\in [a, t),
\\
f(t-0), & x=t,
\\
f(x),           &x\in (t, b],
\end{array}
\right.
\\
f_{v}(x):\!\!\!\! &=& \!\!\!\! \left\{
\begin{array}{ll}
f(t+0), & x=t,
\\
f(x),           &x\in (t, b],
\end{array}
\right.
\\
g_{e}(x):\!\!\!\! &=& \!\!\!\! \left\{
\begin{array}{ll}
g(x),           &x\in [c,s),
\\
g(s-0), & x=s,
\\
g(x),           &x\in (s,d],
\end{array}
\right.
\\
g_{v}(x):\!\!\!\! &=& \!\!\!\! \left\{
\begin{array}{ll}
g(s+0), & x=s,
\\
g(x).           &x\in (s,d].
\end{array}
\right.
\end{eqnarray*}
Then define sequence $\{f_{e}^{n}(t)\}$, where $n=0,1,2,\cdots$. Obviously, $f_{e}^{2k}(t)\rightarrow a$, and $f_{e}^{2k+1}(t)\rightarrow b$. With these sequence, give a partition of $I$, i.e.,
\begin{eqnarray*}
(a,b)=\cup_{k=0}^{+\infty}[f_{e}^{2k+2}(t),f_{e}^{2k}(t)]\cup[t,f_{e}(t)]\cup\cup_{k=0}^{+\infty}[f_{e}^{2k+1}(t),f_{e}^{2k+3}(t)].
\end{eqnarray*}
Choose any increasing homeomorphism $\psi_{0}:[f_{e}^{2}(t),t]\rightarrow [g_{e}^{2}(s),s]$ as an initial function such that $\psi_{0}$ goes through points $(f(t),g(s))$ and $(f_{v}(t),g_{v}(s))$.Ddefine for $n=1,~2,~3, \cdots$
\begin{eqnarray*}
\psi_{1}(x):\!\!\!\! &=& \!\!\!\! \left\{
\begin{array}{ll}
c, &x=a,
\\
\psi_{0}(x),           &x\in (f_{e}^{2}(t),t],
\\
g_{e}^{n}\circ\psi_{0}\circ f_{e}^{-n}(x), &x\in (f_{e}^{n+2}(t),f_{e}^{n}(t)], or [f_{e}^{n}(t),f_{e}^{n+2}(t)),
\\
d, &x=b,
\end{array}
\right.
\end{eqnarray*}
\begin{eqnarray*}
\psi_{2}(x):\!\!\!\! &=& \!\!\!\!
g_{v}^{-1}\circ\psi_{0}\circ f_{v}(x),~~x\in (t,f_{e}(t)),
\end{eqnarray*}
and
\begin{eqnarray*}
\psi(x):\!\!\!\! &=& \!\!\!\! \left\{
\begin{array}{ll}
\psi_{1}(x),           &x\in [a,t)\cup[f_{e}(t),b],
\\
s,           &x=t,
\\
\psi_{2}(x),           &x\in (t,f_{e}(t)).
\end{array}
\right.
\end{eqnarray*}
One can see that $\psi: [a,b]\to [c,d]$ is an increasing conjugacy from $f$ to $g$, and every conjugacy can be obtained in this manner.

For subcase \textbf{(ii)}, define functions $\tilde{f}_{e}$, $\tilde{f}_{v}$, $\tilde{g}_{e}$ and $\tilde{g}_{v}$ as follows
\begin{eqnarray*}
\tilde{f}_{e}(x):\!\!\!\! &=& \!\!\!\! \left\{
\begin{array}{ll}
f(x),           &x\in [a, t),
\\
f(t-0), & x=t,
\\
f(x),           &x\in (t, b],
\end{array}
\right.
\\
\tilde{f}_{v}(x):\!\!\!\! &=& \!\!\!\! \left\{
\begin{array}{ll}
f(t+0), & x=t,
\\
f(x),           &x\in (t, b],
\end{array}
\right.
\\
\tilde{g}_{e}(x):\!\!\!\! &=& \!\!\!\! \left\{
\begin{array}{ll}
g(x),           &x\in [c,s),
\\
g(s-0), & x=s,
\end{array}
\right.
\\
\tilde{g}_{v}(x):\!\!\!\! &=& \!\!\!\! \left\{
\begin{array}{ll}
g(x), &x\in [c,s),
\\
g(s+0),  &x=s,
\\
g(x).           &x\in (s,d].
\end{array}
\right.
\end{eqnarray*}
 Then define sequence $\{\tilde{f}_{e}^{n}(t)\}$, where $n=0,1,2,\cdots$. Obviously, $\tilde{f}_{e}^{2k}(t)\rightarrow a$, and $\tilde{f}_{e}^{2k+1}(t)\rightarrow b$. With these sequence, give a partition of $I$, i.e.,
\begin{eqnarray*}
(a,b)=\cup_{k=0}^{+\infty}[\tilde{f}_{e}^{2k+2}(t),\tilde{f}_{e}^{2k}(t)]\cup[t,\tilde{f}_{e}(t)]\cup\cup_{k=0}^{+\infty}[\tilde{f}_{e}^{2k+1}(t),\tilde{f}_{e}^{2k+3}(t)].
\end{eqnarray*}
 Choose any decreasing homeomorphism $\tilde{\psi}_{0}:[\tilde{f}_{e}^{2}(t),t]\rightarrow [s,\tilde{g}_{v}^{2}(s)]$ as an initial function such that $\tilde{\psi}_{0}$ goes through points $(f(t),g(s))$ and $(\tilde{f}_{v}(t),\tilde{g}_{e}(s))$.
Define for $n=1,~2,~3, \cdots$
\begin{eqnarray*}
\tilde{\psi}_{1}(x):\!\!\!\! &=& \!\!\!\! \left\{
\begin{array}{ll}
d, &x=a,
\\
\tilde{\psi}_{0}(x),           &x\in (\tilde{f}_{e}^{2}(t),t],
\\
\tilde{g}_{v}^{n}\circ\tilde{\psi}_{0}\circ \tilde{f}_{e}^{-n}(x), &x\in (\tilde{f}_{e}^{n+2}(t),\tilde{f}_{e}^{n}(t)], or [\tilde{f}_{e}^{n}(t),\tilde{f}_{e}^{n+2}(t)),
\\
c, &x=b,
\end{array}
\right.
\end{eqnarray*}
\begin{eqnarray*}
\tilde{\psi}_{2}(x):\!\!\!\! &=& \!\!\!\!
\tilde{g}_{e}^{-1}\circ\tilde{\psi}_{0}\circ \tilde{f}_{v}(x),~~x\in (t,\tilde{f}_{e}(t)),
\end{eqnarray*}
and
\begin{eqnarray*}
\tilde{\psi}(x):\!\!\!\! &=& \!\!\!\! \left\{
\begin{array}{ll}
\tilde{\psi}_{1}(x),           &x\in [a,t)\cup[\tilde{f}_{e}(t),b],
\\
s,           &x=t,
\\
\tilde{\psi}_{2}(x),           &x\in (t,\tilde{f}_{e}(t)).
\end{array}
\right.
\end{eqnarray*}
We have $\tilde{\psi}: [a,b]\to [c,d]$ is a decreasing conjugacy from $f$ to $g$ and every conjugacy can be obtained in this manner.
 $\Box$

\section{Smoothness of conjugacy}
\setcounter{equation}{0}
\setcounter{lm}{0}
\setcounter{thm}{0}
\setcounter{rmk}{0}
\setcounter{df}{0}
\setcounter{cor}{0}
\setcounter{exa}{0}

For convenience, we only consider subcase \textbf{(i)} under case \textbf{(C1)} in Theorem \ref{increasing-case-l}, Others are similar.
\begin{thm}
The conjugacy $\varphi$ in the proof of Theorem \ref{increasing-case-l} is continuously differentiable, provided the
following hypotheses are added:
\begin{description}

\item[(a)]
$f_l$, $f_r$, $g_l$ and $g_r$ are continuously differentiable,
$f_{l}^{\prime}(x)\neq 0$ for $x\in[a,t]$ and $f_{r}^{\prime}(x)\neq 0$ for $x\in[t,b]$;
$\varphi_{l_{0}}$ and $\varphi_{r_{0}}$ are continuously differentiable on $[f_{l}(t),t]$ and $[t,f_{r}(t)]$, respectively, satisfying $\varphi_{l_{0}}^{\prime}(t)=\varphi_{r_{0}}^{\prime}(t)$ and
$$\varphi^{\prime}_{l_{0}}(f_l(t))=\frac{g_l^{\prime}(s)}{f_l^{\prime}(t)}\varphi_{l_{0}}^{\prime}(t),$$
$$\varphi_{r_{0}}^{\prime}(f_{r}(t))=\frac{g_r^{\prime}(s)}{f_r^{\prime}(t)}\varphi_{r_{0}}^{\prime}(t);$$

\item[(b)] For any $x\in(f_l(t),t]$, $\varphi_{l_0}^{\prime}(x)\prod_{j=0}^{\infty}\frac{g_{l}^{\prime}(g_{l}^{j}(\varphi_{l_0}(x)))}{f_{l}^{\prime}(f_{l}^{j}(x))}=L_{1}$;

\item[(c)]
For any $x\in[t,f_r(t))$, $\varphi_{r_0}^{\prime}(x)\prod_{j=0}^{\infty}\frac{g_{r}^{\prime}(g_{r}^{j}(\varphi_{r_0}(x)))}{f_{r}^{\prime}(f_{r}^{j}(x))}=L_{2}$.

\end{description}

\label{smoothness-conjugacy}
\end{thm}
\noindent{\bf Proof.}
By induction, we can prove $\varphi_{l}$ is continuously differentiable on the interval $(a,t]$.
Now it suffices to show $\varphi_{l}$ is continuously differentiable at $x=a$. Since $\varphi_{l}$ is continuously differentiable on the interval $(a,~t]$, we have
$$
\varphi_{l}^{\prime}(f_{l}(x))f_{l}^{\prime}(x)=g_{l}^{\prime}(\varphi_{l}(x))\varphi_{l}^{\prime}(x), x\in (a,~t].
$$
It follows from $f_{l}^{\prime}(x)\neq0$ that
\begin{eqnarray*}
\varphi_{l}^{\prime}(f_{l}(x))=\frac{g_{l}^{\prime}(\varphi_{l}(x))}{f_{l}^{\prime}(x)}\varphi_{l}^{\prime}(x), x\in (a,~t].
\label{iterative2}
\end{eqnarray*}
By induction, we have
\begin{eqnarray}
\varphi_{l}^{\prime}(f_{l}^{n}(x))=\varphi_{l}^{\prime}(x)\prod_{j=1}^{n-1}\frac{g_{l}^{\prime}(g_{l}^{j}(\varphi_{l}(x)))}{f_{l}^{\prime}(f_{l}^{j}(x))}, x\in (a,~t].
\end{eqnarray}
In particular, we choose $x\in (f_l(t),t]$ and have
\begin{eqnarray}
\varphi_{l}^{\prime}(f_{l}^{n}(x))=\varphi_{l_{0}}^{\prime}(x)\prod_{j=0}^{n-1}\frac{g_{l}^{\prime}(g_{l}^{j}(\varphi_{l_{0}}(x)))}{f_{l}^{\prime}(f_{l}^{j}(x))}, \ \ x\in (f_l(t),t].
\label{iterative3}
\end{eqnarray}
By the assumption {\bf (b)}, one has
$\varphi_{l_{0}}^{\prime}(x)\prod_{j=0}^{n-1}\frac{g_{l}^{\prime}(g_{l}^{j}(\varphi_{l_{0}}(x)))}{f_{l}^{\prime}(f_{l}^{j}(x))}\rightarrow L_1 $ as $n\rightarrow\infty$ for $x\in (f_l(t),t]$. Then
\begin{eqnarray}
\lim_{n\rightarrow\infty}\varphi_{l}^{\prime}(f_{l}^{n}(x))=L_1 ,\quad  x\in (f_l(t),t].
\label{limit112}
\end{eqnarray}
By the Mean Value Theorem, for any $x\in(a,t]$, there exists $\xi\in (a,x)$ such that
\begin{align*}
\frac{\varphi_{l}(x)-\varphi_{l}(a)}{x-a}=\varphi_{l}^{\prime}(\xi).
\end{align*}
It follows from (\ref{limit112}) that
$$
\varphi_{l}^{\prime}(a)=\lim_{x\rightarrow a^+}\frac{\varphi_{l}(x)-\varphi_{l}(a)}{x-a}=\lim_{\xi\rightarrow a^+}\varphi_{l}^{\prime}(\xi)=L_1.
$$
Therefore, $\varphi_{l}$ is continuously differentiable on $[a,t])$.
Similarly, we can prove $\varphi_{r}$ is also continuously differentiable on $[t,b]$ by the condition \textbf{(c)}.
Thus $\varphi$ is a continuously differentiable conjugacy. $\Box$

\begin{rmk}
 We can also investigate the smoothness of conjugacies in Theorem \ref{decreasing-case-1} as Theorem \ref{smoothness-conjugacy}. More hypotheses are imposed on the initial function $\psi_{0}$.
 \end{rmk}

\section{Applications to Lorenz maps}
\begin{figure}[H]
\begin{minipage}[t]{0.5\linewidth}
\centering
\includegraphics[width=0.8\textwidth]{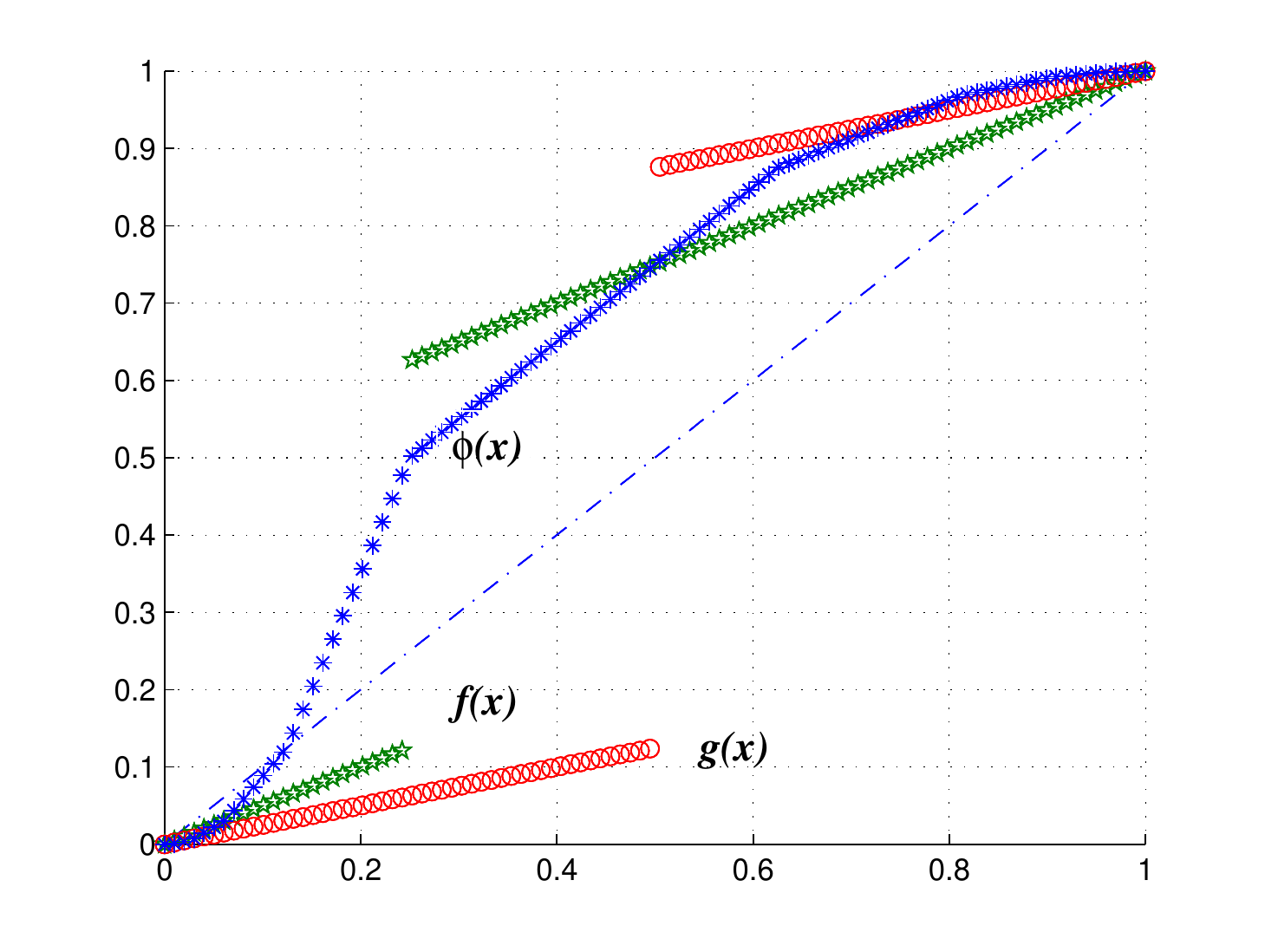}
\caption{an increasing case}
\label{exm}
\end{minipage}%
\begin{minipage}[t]{0.5\linewidth}
\centering
\includegraphics[width=1.1\textwidth]{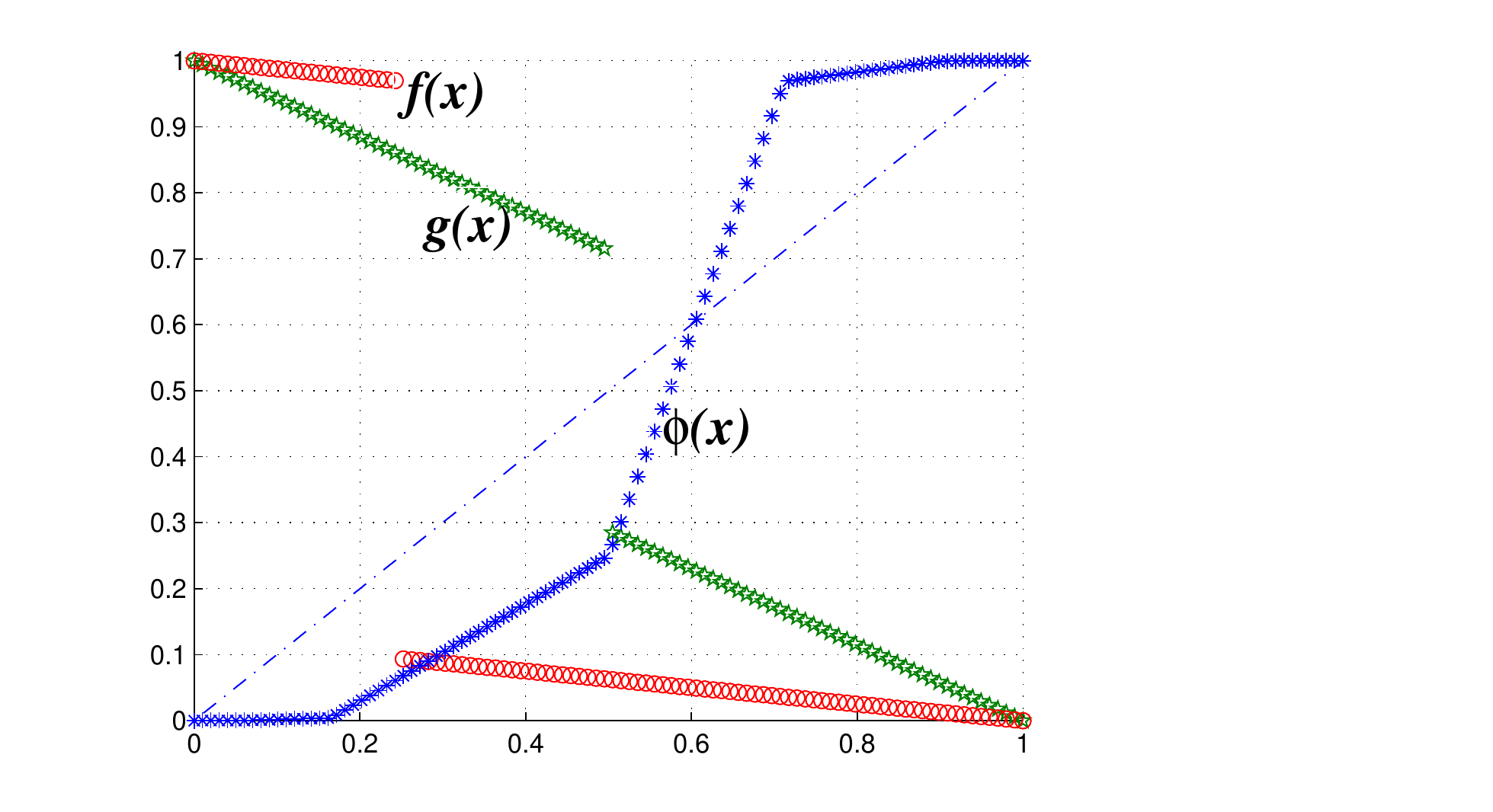}
\caption{a decreasing case}
\label{hfg}  
\end{minipage}
\end{figure}

\begin{exa}
Consider the following two Lorenz maps
\begin{eqnarray*}
f(x)=\left\{
\begin{array}{ll}
\frac{1}{2}x,           &x\in [0,\frac{1}{4}),
\\[1mm]
\frac{3}{16},           &x=\frac{1}{4},
\\[1mm]
\frac{1}{2}x+\frac{1}{2}, &x\in (\frac{1}{4},1],
\end{array}
\right.
\quad\qquad
g(x)= \left\{
\begin{array}{ll}
\frac{1}{4}x,           &x\in [0,\frac{1}{2}),
\\[1mm]
\frac{5}{16},           &x=\frac{1}{2},
\\[1mm]
\frac{1}{4}x+\frac{3}{4}, &x\in (\frac{1}{2},1].
\end{array}
\right.
\end{eqnarray*}
One can see that  $f\in\mathcal{A}_{\frac{1}{4}}([0,1])$ and $g\in\mathcal{A}_{\frac{1}{2}}([0,1])$.
By the proof of Theorem \ref{increasing-case-l}, a topological conjugacy $\phi$ from $f$ to $g$ is constructed in Fig.\ref{exm}.
\end{exa}
\begin{exa}
Consider the following two Lorenz maps
\begin{eqnarray*}
f(x):\!\!\!\! &=& \!\!\!\! \left\{
\begin{array}{ll}
-\frac{23}{40}x+1,           &x\in [0,\frac{1}{2}),
\\
\frac{17}{40},              &x=\frac{1}{2},
\\
-\frac{23}{40}x+\frac{23}{40}, &x\in (\frac{1}{2},1],
\end{array}
\right.
\quad\qquad
g(x)=\left\{
\begin{array}{ll}
-\frac{1}{8}x+1,           &x\in [0,\frac{1}{4}),
\\[1mm]
\frac{53}{272},             &x=\frac{1}{4},
\\[1mm]
-\frac{1}{8}x+\frac{1}{8}, &x\in (\frac{1}{4},1],
\end{array}
\right.
\end{eqnarray*}
It is clear that $f\in\mathcal{B}_{\frac{1}{2}}([0,1])$ and $g\in\mathcal{B}_{\frac{1}{4}}([0,1])$.
By the proof of Theorem \ref{decreasing-case-1}, a topological conjugacy $\phi$ from $f$ to $g$ is constructed in Fig.\ref{hfg}.
\end{exa}

\end{document}